\documentclass[12pt]{article}
\usepackage{amssymb,amsmath,amscd,amsthm}
\usepackage{graphicx,epsfig, color,float}

\usepackage{graphicx}
\usepackage[active]{srcltx}

\newtheorem{theorem}{Theorem}[section]

\newtheorem{lemma}[theorem]{Lemma}

\def\be{\color{black}}

\date{}

\begin{document}

\date{}
\title{Asymptotics in the Dirichlet Problem for Second Order Elliptic Equations with Degeneration on the Boundary}
\author{
 M. Freidlin\footnote{Dept of Mathematics, University of Maryland,
College Park, MD 20742, mif@math.umd.edu}, L.
Koralov\footnote{Dept of Mathematics, University of Maryland,
College Park, MD 20742, koralov@math.umd.edu}
%\footnote{Dept of Mathematics, University of Maryland, College
%Park, MD 20742, koralov@math.umd.edu}
%\footnote{Dept of Mathematics, University of Maryland,
%College Park, MD 20742, mif@math.umd.edu}
} \maketitle

\begin{abstract}  
We study small perturbations of the Dirichlet problems for second order elliptic equations that degenerate on the boundary.
The limit of the  solution, as the perturbation tends to zero, is calculated.  The result is based on a certain asymptotic
self-similarity near the boundary, 
which holds in the generic case. In the last section, we briefly consider the stabilization of solutions to the
corresponding parabolic equations with a small parameter. Metastability effects arise in this case: the asymptotics of the solution depends
on the time scale. Initial-boundary value problem  with the Neumann boundary condition is discussed in the last section as well. 
\end{abstract}



\section{Introduction} \label{intro}

Let $D$ be a bounded domain in $\mathbb{R}^d$ with a sufficiently smooth boundary $\partial D = S$. In most of the paper,
we assume that $S$ is connected; the case when $S$ consists of a finite number of connected components is discussed in Section~\ref{gener}.

Let $v_0,...,v_d$ be sufficiently smooth vector fields on $\mathbb{R}^d$. Define the operator
\begin{equation} \label{opl}
Lu (x) = L_0 u (x) + \frac{1}{2} \sum_{i=1}^d L_i^2 u (x),~~~x \in \mathbb{R}^d,
\end{equation}
where $L_i$ is the operator of differentiation along the vector field $v_i$, $i =0,...,d$. The operator is the generator of 
the diffusion process  $X^x_t$, $x \in \mathbb{R}^d$, defined by the stochastic differential equation
\[
d X^{x}_t = v_0( X^{x}_t) dt + \sum_{i=1}^d v_i( X^{x}_t) \circ d W^i_t,~~~~~X^{x}_0 = x,
\]
where $W^i_t$ are independent Wiener processes, and the stochastic term is understood in the Stratonovich sense.  The Stratonovich form for the SDE
and the corresponding way to write the differential operator are convenient here since they allow one to provide a coordinate-independent description of the process. 

For each $x \in S$, we define $\mathbf{n}(x)$ to be the unit interior normal vector to $S$ at $x$, and define $T(x)$ to be the tangent space to $S$ at $x$.  We assume that: 

(a) ${\rm span} (v_0(x),v_1(x),...,v_d(x)) = {\rm span} (v_1(x),...,v_d(x)) = T(x)$ for $x \in S$;

(b) ${\rm span} (v_1(x),...,v_d(x)) = \mathbb{R}^d$ for $x \in D$. 
\\
Assumption (a)~implies that the surface $S$ is invariant for the process, and, consequently, the operator $L$ can be applied to functions
defined on $S$. Moreover, when restricted to the surface, $L$ is uniformly elliptic, and so the process is ergodic on $S$. We denote the
unique invariant probability measure on $S$ by $\pi$. Assumption (b) implies that $L$ is uniformly elliptic on each compact subdomain of $D$.

Consider now a small non-degenerate perturbation $L^\varepsilon$ of the operator $L$:
\[
L^\varepsilon = L + \varepsilon^2  \tilde{L},~~~\varepsilon > 0,
\] 
with
\[
\tilde{L}u = \tilde{L}_0 + \frac{1}{2} \sum_{i=1}^d \tilde{L}_i^2.
\]
Here, $\tilde{L}_i$ is the operator of differentiation along the vector field $\tilde{v}_i$, $i =0,...,d$.  
In order to make our assumption on the non-degeneracy of the
perturbation more precise, we state it as follows: 

(c) ${\rm span} (\tilde{v}_1(x),...,\tilde{v}_d(x)) = \mathbb{R}^d$ for $x \in \mathbb{T}^d$. 
\\
The operator $L^\varepsilon$ is the generator for the   diffusion process  $X^{x,\varepsilon}_t$ that satisfies
\[
d X^{x,\varepsilon}_t = (v_0 + \varepsilon^2   \tilde{v}_0) (X^{x,\varepsilon}_t ) dt + \sum_{i=1}^d v_i(X^{x,\varepsilon}_t ) \circ d W^i_t + \varepsilon 
\sum_{i =1}^d \tilde{v}_i (X^{x,\varepsilon}_t ) \circ d \tilde{W}^i_t,~~~~~X^{x,\varepsilon}_t = x,
\]
where $\tilde{W}^i_t$ are independent Wiener processes (also independent of all $W^i_t$).

Since the operator $L^\varepsilon$ is uniformly elliptic on $\overline{D}$, and $S  = \partial D$ is smooth, the Dirichlet problem
\begin{equation} \label{direq}
L^\varepsilon u^\varepsilon (x)  = 0,~~x \in D;~~u^\varepsilon(x)  = \psi(x),~~x \in \partial D,
\end{equation}
has a unique solution for each continuous boundary function $\psi$ for each $\varepsilon > 0$. The main result of the current paper
concerns the behavior of the solution $u^\varepsilon$ as $\varepsilon \downarrow 0$. Note that the solution admits the probabilistic 
representation
\[
u^\varepsilon (x) = \mathrm{E} \psi (X^{x,\varepsilon}_{\tau^{x,\varepsilon}(\partial D)}),
\]
where $\tau^{x,\varepsilon}(\partial D)$ is the first time when $X^{x,\varepsilon}_t$ reaches $\partial D$. Therefore, the asymptotic behavior of $u^\varepsilon$ can be understood by studying the limiting distribution of 
$X^{x,\varepsilon}_{\tau^{x,\varepsilon}(\partial D)}$ as $\varepsilon \downarrow 0$.

A probabilistic representation can also be used for the solutions of the Neumann and the (parabolic) initial-boundary value problems 
for the operator $L^\varepsilon$. 
When $S$ has multiple connected components, the asymptotic behavior of the solutions to the Neumann and 
initial-boundary value problems is closely related to the question about the metastable distributions of 
the process $X^{x,\varepsilon}_t$ (or the corresponding process with reflection on the boundary in the case of the Neumann problem), i.e., about the limiting behavior of the process as $\varepsilon \downarrow 0$ and,
simultaneously, $t = t(\varepsilon) \rightarrow \infty$. In certain cases, the limiting distribution for $X^{x,\varepsilon}_{t(\varepsilon)}$
depends on the way in which $t(\varepsilon)$ approaches infinity. Such problems are briefly discussed in Section~\ref{gener} and will be the subject of a forthcoming paper.

Let us breifly mention the connection between equation (\ref{direq}) studied in our paper and equations without the regularizing 
term $\varepsilon^2 \tilde{L}$ added to the operator. Elliptic and parabolic equations with coefficients 
degenerating on the boundary and the corresponding diffusion processes were first analyzed in one-dimensional case by Feller (\cite{Fel},
\cite{Fel2}). The results were generalized to multi-dimensional case by Hasminskii (\cite{Has}), who provided certain necessary and 
(separately) sufficient conditions, in terms of the
coefficients of the generator, for the boundary to be attracting or repelling for the corresponding diffusion process. The notion of
repelling or attracting boundary also plays an important role in our situation and does not depend on the regularizing operator $\tilde{L}$. 
However, in our situation, the diffusion along the boundary remains  non-degenerate, and the asymptotics of the elliptic problem is
determined by the interplay between  this diffusion and the perturbation. Even if the boundary is repelling for the unperturbed
diffusion, the presence of the perturbation means that the process governed by $L^\varepsilon$ reaches the boundary, and its behavior near
the boundary needs to be studied in oreder to understand the asymptotics of the solution.

Finally, let us mention that elliptic and parabolic equations with degeneration inside the domain were studied under various assumptions on
the coefficients (see \cite{F85}, \cite{OR} and references therein). The first results on small perturbations of such equations  concerned the case when the unperturbed operator was of the first order (\cite{Lev}, \cite{H2}, \cite{F85}, \cite{FW}).

The paper is organized as follows. In Section \ref{anr}, we state the assumptions in more detail and formulate the main result. The limit of $u^\varepsilon$ is expressed in terms of  the solution to an auxiliary elliptic problem, with the operator obtained from $L^\varepsilon$ by using the leading terms in the asymptotic expansion of the coefficients in the vicinity of $S$:  after a change of variables, the operator and the corresponding diffusion process can be considered on the space $S \times [0, \infty)$, where the 
second coordinate corresponds to the distance (scaled by $\varepsilon$) of a point to the  surface $S$. The resulting operator is 
homogeneous in the second variable. 
Such problems and the corresponding diffusion processes are studied in  Section \ref{homog}. The main
result is proved in Section~\ref{mresp}. Generalizations concerning the case of a domain whose bounary has multiple connected components are discussed  in  Section~\ref{gener}. 

\section{Structure of the operators $L$ and $L^\varepsilon$ near the boundary. Formulation of the main result} \label{anr}

We assume that the boundary $S$ is $C^4$-smooth, $v_0,...,v_n, \tilde{v}_0,...,\tilde{v}_n \in C^3(\bar{D})$, and
that (a)-(c) hold. Let us specify how the coefficients degenerate on $S$. 
Roughly speaking, while there is no diffusion or drift across the boundary, the diffusion should degenerate in a generic way, i.e., the normal derivative of the diffusion coefficient in the direction orthogonal to $S$ should be non-zero. Let us make this assumption more precise.  

Recall that $\mathbf{n}$ is the field of unit interior normal vectors to $S$. For each $x \in \bar{D}$ in a sufficiently small neighborhood of $S$, there is a unique $y(x) \in S$ such 
that ${\rm dist}(x, y(x)) = {\rm dist}(x, S)$. Define
\[
z(x) = \langle x -  y(x), \mathbf{n}(y(x)) \rangle.
\]
Thus $z$ measures the distance of $x$ from $S$. For a sufficiently 
small $\delta > 0$, this function is defined
and belongs to $C^3(S^\delta)$, where $S^\delta$ is the $\delta$-neighborhood of $S$ in $\bar{D}$. Observe that 
$\varphi(x) = (y(x), z(x))$ is a bijection between $S^\delta$ and the set $S \times [0, \delta)$, i.e., $(y,z)$ can be viewed as a new set of coordinates on  $S^\delta$. 

From the conditions placed on the vector fields 
$v_0,v_1,...,v_d$, it follows that the operator $L$ can be applied to functions defined on $S$. 
In order to stress that we are considering the restriction of the operator to $S$ (where variables $y$ are used), we denote the resulting operator by $L_y$. 
Let us write the processes $X^x_t$ and $X^{x,\varepsilon}_t$ in $(y,z)$ coordinates. Note that the operator
$L_y$ acting in the $y$ variables can be applied to functions of $(y,z)$ by treating $z$ as a parameter. 
\begin{lemma} \label{gene} The generator of the process $X^x_t$ in $(y,z)$ coordinates can be written as:
\begin{equation} \label{opr}
L u  = L_y u  +  z^2 \alpha(y) \frac{\partial^2 u}{\partial z^2} +   z  \beta(y) \frac{\partial u}{\partial z}  + z {\mathcal{D}}_y \frac{\partial u}{\partial z} + R u 
\end{equation}
with
\begin{equation} \label{cott}
R u = z {\mathcal{K}}_y u   + z^2 \mathcal{N}_y \frac{\partial u}{\partial z}  +  z^3 \sigma(y,z) \frac{\partial^2 u}{\partial z^2},
\end{equation}
where ${\mathcal{D}}_y$ is a  differential operator with first-order derivatives in $y$, whose coefficients depend only on the $y$ variables, ${\mathcal{K}}_y$ is a  differential operator  on $S \times [0, \delta)$ with first- and second-order derivatives in $y$,  ${\mathcal{N}}_y$ is a differential operator on 
$S \times [0, \delta)$ with first-order derivatives in $y$ and a potential term, All the 
operators have continuously differentiable coefficients, 
while $\alpha, \beta \in C^1(S)$,  $\sigma \in C^1 (S^\delta)$. 
The generator of the process $X^{x,\varepsilon}_t$ in $(y,z)$ coordinates is the operator
\[
L^\varepsilon = L + \varepsilon^2   \tilde{L},
\]
where $\tilde{L}$ is a second-order uniformly elliptic differential operator with continuously differentiable coefficients in $(y,z)$ variables. 
\end{lemma}
\proof
Consider the simplest case: $S$ is one-dimensional (i.e., $d = 2$) and defined, in the $(y,z)$ coordinates, in a 
neighborhood of a point $(y_0, z_0) \in S$, as the line $\{(y,z): z = 0\}$.  (The case when $d > 2$ requires only slightly more complicated notations, while
the case when $S$ is a curve (surface) can be reduced to the case when it is a linear subspace by a change of variables.)  Each term of the operator $L$ defined in (\ref{opl}) can
be considered separately, so we can assume that
\[
Lu = \frac{1}{2} \frac{\partial}{\partial v}(\frac{\partial u}{\partial v}),
\] 
where $v = (v^1, v^2)$ is a vector field tangent to $S$ (the term with the first order derivative can be considered similarly). Thus
\[
Lu = \frac{1}{2}(v^1)^2 \frac{\partial^2 u}{\partial y^2} +  v^1 v^2  \frac{\partial^2 u}{\partial y \partial z} + 
\frac{1}{2} (v^2)^2 \frac{\partial^2 u}{\partial z^2} + 
\frac{1}{2} (v^1 \frac{\partial v^1}{\partial y} + v^2 \frac{\partial v^1}{\partial z}) \frac{\partial u}{\partial y} + \frac{1}{2} (v^1 \frac{\partial v^2}{\partial y} +
v^2 \frac{\partial v^2}{\partial z}) \frac{\partial u}{\partial z}.
\]
Using smoothness of $v^1, v^2$ and the fact that $v^2(y,0) = 0$, we can write
\[
v^1(y,z) = v^1(y,0) +  g_1(y,z)z,~~v^2(y,z) = \frac{\partial v^2}{\partial z} (y,0) z +  g_2(y,z)z^2,~~
\]
\[
\frac{\partial v^1}{\partial y}(y,z) = \frac{\partial v^1}{\partial y}(y,0) +  g_3(y,z)z,~~ 
\frac{\partial v^1}{\partial z}(y,z) = \frac{\partial v^1}{\partial z}(y,0) +  g_4(y,z)z,
\]
\[
\frac{\partial v^2}{\partial y}(y,z) = \frac{\partial^2 v^2}{\partial y \partial z} (y,0) z +  g_5(y,z)z^2,~~
\frac{\partial v^2}{\partial z}(y,z) = 
\frac{\partial v^2}{\partial z}(y,0) +  g_6(y,z)z,
\]
where $g_1,...,g_6$ are smooth functions. Expressing the coefficients of $L$ using these expansions, we obtain the desired form of the operator. 
The statement about the form of $L^\varepsilon$ follows immediately.
\qed
%\[
% \frac{1}{2}(v^1)^2(y,0) \frac{\partial^2 u}{\partial y^2} +  v^1 (y,0) z \frac{\partial v^2}{\partial z} (y,0)  \frac{\partial}{\partial y}( 
%\frac{\partial u}{ \partial z})  
%\]
%\[
%\frac{1}{2}\left( (v^1)^2(y,z)- (v^1)^2(y,0) \right) \frac{\partial^2 u}{\partial y^2} +  \left(v^1 v^2(y,z) - v^1 (y,0)  \frac{\partial v^2}{\partial z} (y,0) z \right)
% \frac{\partial}{\partial y}( 
%\frac{\partial u}{ \partial z})  
%\]
\\

\noindent
{\bf Remark.} Define $h_1 =  L z$, $h_2 = \frac{1}{2} Lz^2$. Then
\[
h_1(x) = \beta(y(x)) z(x) + O(z^2(x)),~~h_2(x) = \alpha(y(x)) z^2(x) + O(z^3(x)),~~{\rm as}~~z(x) \downarrow 0, 
\]
which provides a simple way to identify $\alpha(y)$ and $\beta(y)$.
\\

The assumption that the diffusion in the direction orthogonal to the boundary degenerates in a generic way is expressed by the 
requirement that $\alpha > 0$ for each $y \in S$ (we could weaken this assumption and instead assume that there is $y \in S$ such that 
$\alpha(y) > 0$).  
Recall that $\pi$ is the unique invariant probability measure on $S$ for the diffusion $X^x_t$, which is non-degenerate when restricted to $S$.
Define 
\[
\bar{\alpha} = \int_{S} \alpha(y) d\pi(y),~~\bar{\beta} = \int_{S} \beta(y) d\pi(y) .
\]
We will see (Lemma~\ref{mlem}) that  if $\bar{\alpha} > \bar{\beta} $,
then $\mathrm{P} (\lim_{t \rightarrow \infty} {\rm dist}(X^x_t, S) = 0) = 1$ for each $x \in \overline{D}$. 
If $\bar{\alpha} < \bar{\beta}$, then this
probability is zero unless $x \in S$. We will
refer to $S$ as attracting if $\bar{\alpha} > \bar{\beta} $, repelling if $\bar{\alpha} < \bar{\beta}$, and neutral if
 $\bar{\alpha} = \bar{\beta} $.
\\

Formula (\ref{opr}) shows that the generator of $X^{x}_t$ can be approximated near the boundary by an operator that is homogeneous in the variable $z$. The  remainder term with the operator $R$ can be made small by considering a sufficiently small neighborhood of $S$. The same approximation is useful for the generator of $X^{x, \varepsilon}_t$, except in a yet smaller neighborhood of $S$, where the 
perturbation $\varepsilon^2   \tilde{L}$ is comparable to or larger than the operator $L$, since the coefficients of the latter degenerate near $S$. In order to understand the behavior of $X^{x, \varepsilon}_t$ in $S^{r \varepsilon}$, 
we introduce the change of variables
\[
\Psi_\varepsilon(y,z) = (y,\frac{z}{\varepsilon}),~~~\Psi_\varepsilon: S \times[0,\delta) \rightarrow S \times [0,\frac{\delta}{\varepsilon})
\subset S \times [0,\infty)
\]
and the operator $M^\varepsilon u  =  {L}^\varepsilon (u (\Psi_\varepsilon)) (\Psi_\varepsilon^{-1})$.
This operator is the generator of the  process 
\begin{equation} \label{chp}
\mathcal{X}^{{\bf x},\varepsilon}_t := \Psi_\varepsilon (X^{\Psi_\varepsilon^{-1} ({\bf x}), \varepsilon}_t)
\end{equation}
 on $S \times [0,{\delta}/{\varepsilon})$.  To stress the difference between the two sets of coordinates, we use the notation
${\bf x} = (y,{\bf z})$ for the new variables instead of $x = (y,z)$. 
By Lemma~\ref{gene},  
\begin{equation} \label{genpp0}
M^\varepsilon u  = L_y u  +  ({\bf z}^2 \alpha(y) + \rho(y))\frac{\partial^2 u}{\partial {\bf z}^2} +   
{\bf z} \beta(y) \frac{\partial u}{\partial {\bf z}}  + 
{\bf z} {\mathcal{D}}_y \frac{\partial u}{\partial {\bf z}} + 
\widehat{R}^\varepsilon u =: M u +  \widehat{R}^\varepsilon u, 
\end{equation}
where $\rho$ is the coefficient at the second derivative in the variable $z$ at $z=0$ in the operator $\tilde{L}$, and $\widehat{R}^\varepsilon$ is a second order operator with continuously differentiable coefficients 
that tend to zero uniformly on $S \times [0,r]$ as $\varepsilon \downarrow 0$ for each $r > 0$. Thus $M^\varepsilon$ is a small perturbation of the operator $M$, which does not depend on $\varepsilon$. 

We can view $M$ as an operator on $S \times [0,\infty)$. The solutions to the equations in the following lemma
are sought in the spaces $\mathcal{C} : = C^2(S \times (0,\infty)) \bigcap C_b(S \times [0,\infty))$ or 
$\mathcal{C}^2 := C^2(S \times [0,\infty))$. The lemma will be proved in 
Section~\ref{homog}. 
\begin{lemma} \label{mleee}
Let $f \in C(S)$.

(a)  If the boundary $S$ is attracting or neutral, then there is a unique solution $u \in \mathcal{C}$ to the equation
\begin{equation} \label{fieq}
M u (y,{\bf z}) = 0,~~y \in S, {\bf z} > 0;~~~~u(y,0) = f(y).
\end{equation}

(b) If the boundary $S$ is repelling, then there is a unique solution $h \in \mathcal{C}^2$ to the equation
\begin{equation} \label{fieq2}
M h (y,{\bf z}) = 0,~~y \in S, {\bf z} > 0;~~~~h(y,0) \equiv 1;~~~\lim_{{\bf z} \rightarrow \infty} \sup_{y \in S} |h(y,{\bf z})|  = 0.
\end{equation}
There is a unique solution  $u \in \mathcal{C}$ to the equation
\begin{equation} \label{fieq3}
M ( h u) (y,{\bf z}) = 0,~~y \in S, {\bf z} > 0;~~~~u(y,0) = f(y).
\end{equation}

(c) In all the cases, there exists a constant $\overline{u}$ such that  $\overline{u}= \lim_{{\bf z} \rightarrow \infty} u (y,{\bf z})$, uniformly in $y \in S$.
\end{lemma}

\noindent
{\bf Remark.} While Lemma~\ref{mleee} is formulated in PDE terms, it has a simple probabilistic interpretation, as will be seen 
in Section~\ref{homog}. 
Let $\mathcal{X}^{\bf x}_t$ be the process on $S \times [0,\infty)$, starting at ${\bf x} = (y,{\bf z})$, with the generator $M$. 
(Note that we don't have a representation similar to (\ref{chp}) for $\mathcal{X}^{\bf x}_t$ since $M$ was defined by discarding the
higher-order terms in the coefficients of $M^\varepsilon$).  
 Let 
${\tau}^{\bf x} = {\tau}^{(y,{\bf z})} = \inf\{t \geq 0: \mathcal{X}^{\bf x}_t \in S \times \{0\}\} $. In the attracting and neutral cases, 
\[
\overline{u} = \lim_{{\bf z} \rightarrow \infty} \mathrm{E} f(\mathcal{X}^{(y,{\bf z})}_{{\tau}^{(y,{\bf z})}}).
\]
In the repelling case, 
$\mathrm{P} ({\tau}^{(y,{\bf z})} < \infty) = h(y,{\bf z}) \rightarrow 0$ 
as ${\bf z} \rightarrow \infty$, and 
\[
\overline{u} = 
\lim_{{\bf z} \rightarrow \infty} \mathrm{E} (f(\mathcal{X}^{(y,{\bf z})}_{{\tau}^{(y,{\bf z})}})|{\tau}^{(y,{\bf z})} < \infty).
\]
\\

Now we are ready to state the main result of the paper, to be proved in Section~\ref{mresp}.
\begin{theorem} \label{mte}
Let $u^\varepsilon$ be the solution to Dirichlet problem~(\ref{direq}). Then, uniformly on any compact subset
of $D$, 
\[
\lim_{\varepsilon \downarrow 0} u^\varepsilon(x) = \overline{u},
\]
where $\overline{u}$ is determined in Lemma~\ref{mleee}. 
\end{theorem}

\section{Approximation of the process near the boundary} \label{homog}

In the previous section (formula (\ref{genpp0})), we saw that the generator of $\mathcal{X}^{{\bf x},\varepsilon}_t $  in $(y,{\bf z})$ coordinates 
was a perturbation of the operator
\[
M u  = L_y u  +  ({\bf z}^2 \alpha(y) + \rho(y))\frac{\partial^2 u}{\partial {\bf z}^2} +   
{\bf z} \beta(y) \frac{\partial u}{\partial {\bf z}}  + 
{\bf z} {\mathcal{D}}_y \frac{\partial u}{\partial {\bf z}}.  
\]
In fact, the perturbation is small in a sufficiently small neighborhood of $S$. 
Let us examine the behavior of the process with the generator $M$ (which differs from a homogeneous operator by the
presence of the extra term $\rho(y){\partial^2 u}/{\partial {\bf z}^2}$). The process with the generator $M$, earlier denoted by
$\mathcal{{X}}^{\bf x}_t$, will also be written as 
$({Y}^{{\bf x}}_t, {Z}^{{\bf x}}_t)$, where ${\bf x}=(y,{\bf z})$ is the initial point; its state space is $S \times [0,\infty)$;  the
process is stopped upon reaching $S \times \{0\}$. 
Recall that   ${\tau}^{\bf x} = \inf\{t \geq 0: {Z}^{\bf x}_t  = 0\}$ for ${\bf x} \in S \times [0, \infty)$.

\begin{lemma} \label{sttm}
If  $\bar{\alpha} \geq \bar{\beta} $ (the boundary is attracting or neutral), then $\mathrm{P}({\tau}^{\bf x} < \infty) = 1$ for each ${\bf x} = (y,{\bf z}) \in S \times [0, \infty)$. 
If $\bar{\alpha} < \bar{\beta} $ (the boundary is repelling), 
then $\lim_{{\bf z} \rightarrow \infty} \mathrm{P}({\tau}^{(y,{\bf z})} < \infty) = 0$ uniformly in $y \in S$. 
\end{lemma}
\proof Let $\Phi(y,{\bf z}) = (y,\ln({\bf z}))$ be the mapping from $S \times (0,\infty)$ to $S \times\mathbb{R}$.
Consider the process $\Phi(\mathcal{{X}}^{\bf x}_t) = (Y^{\bf x}_t, \ln(Z^{\bf x}_t)) $ 
on $S \times \mathbb{R}$. This process may go to $-\infty$ along the ${\bf z}$-axis in finite time, but this will not cause any problems. \be The generator of this process is 
\begin{equation} \label{genaa}
\mathcal{A} u   = L_y u  +  (\alpha(y) + \rho(y) e^{-2 {\bf z}})(\frac{\partial^2 u}{\partial  {\bf z}^2} - 
\frac{\partial u}{\partial {\bf z}} ) +
\beta(y) \frac{\partial u}{\partial  {\bf z}}  + 
{\mathcal{D}}_y \frac{\partial u}{\partial  {\bf z}}. 
\end{equation} \be
Let $\psi: S \rightarrow \mathbb{R}$ solve
\[
L_y \psi = \alpha - \beta - (\bar{\alpha} - \bar{\beta}),~~~\int_S \psi d \pi  = 0.
\]
Such a function exists and is determined uniquely since $\int_S  (\alpha - \beta - (\bar{\alpha} - \bar{\beta})) d \pi = 0$ and
$\pi$ is the invariant measure for the process with the generator $L_y$. Let $g(y,{\bf z}) =  \psi(y) + \ln({\bf z})$. From the Ito formula applied to
$g(Y^{\bf x}_t, \ln(Z^{\bf x}_t))$, it follows that
\[
h^{{\bf x}}_t : = \psi(Y^{{\bf x}}_t) + 
\ln(Z^{\bf x}_t) - \int_0^t \mathcal{A} g(Y^{{\bf x}}_s, \ln(Z^{{\bf x}}_s)) ds = 
\]
\[
 \psi(Y^{{\bf x}}_t) + 
\ln(Z^{\bf x}_t) +  \int_0^t \rho(Y^{\bf x}_s) (Z^{{\bf x}}_s)^{-2} ds + (\bar{\alpha} - \bar{\beta}) t
\]
is a martingale. For $n \in \mathbb{Z}$, define the following subsets of $ S \times (0,\infty)$:
\begin{equation} \label{subsets}
\Gamma_n  = \{{\bf x}: g({\bf x}) = n\},~~~~ R_n^- = \{{\bf x}: g({\bf x}) \leq n\},~~~~ R_n^+ = \{{\bf x}: g({\bf x}) \geq n\}.
\end{equation}
Since $({Y}^{\bf x}_t, {Z}^{\bf x}_t)$ is a non-degenerate diffusion, there is a constant $c > 0$ such that
\begin{equation} \label{recc}
\mathrm{P} ({\tau}^{\bf x} < \infty) > c,~~{\bf x} \in  R_0^-.
\end{equation}
 For $\bar{\alpha} \geq \bar{\beta} $, 
the process  $\psi(Y^{{\bf x}}_t) + 
\ln(Z^{\bf x}_t) $ is a supermartingale. Since it is unbounded with probability one, the
process $\mathcal{{X}}^{\bf x}_t = (Y^{\bf x}_t, Z^{{\bf x}}_t)$ reaches $R_0^-$ with probability one for each initial point ${\bf x}$. 
From (\ref{recc}) and the strong Markov property, it follows that $\mathrm{P} ({\tau}^{\bf x} < \infty) = 1$ for each ${\bf x}$.  

Now assume that $\bar{\alpha} < \bar{\beta} $. For ${\bf x} \in \Gamma_n$, let 
\[
\sigma^{{\bf x}} = \inf\{t \geq 0: \mathcal{{X}}^{\bf x}_t \in \Gamma_{n-1} \bigcup \Gamma_{n+1}\}.
\]
It is sufficient to show that there exists $c > 0$ such that, for all sufficiently large~$n$,
\begin{equation} \label{slnn}
\mathrm{P} ( \mathcal{{X}}^{\bf x}_{\sigma^{{\bf x}} }  \in \Gamma_{n+1}) \geq \frac{1}{2} + c,~~~
{\bf x} \in \Gamma_n.
\end{equation}
Since $(Y^{\bf x}_t, \ln(Z^{\bf x}_t))$ is a diffusion with coefficients that are bounded on $R^+_0$, 
there is $c' > 0$ such
that, 
\[
\mathrm{E} \sigma^{{\bf x}} \geq c',~~~
{\bf x} \in \Gamma_n,
\]
provided that $n \geq 1$. Since $\mathrm{E}(h^{{\bf x}}_{\sigma^{{\bf x}}} - h^{{\bf x}}_0) = 0$,
\begin{equation} \label{innn}
2 \mathrm{P} ( \mathcal{{X}}^{\bf x}_{\sigma^{\bf x}}  \in \Gamma_{n+1}) - 1 
- \mathrm{E} \int_0^{ \sigma^{{\bf x}}} \left( (\bar{\beta} - \bar{\alpha} ) - \rho(Y^{\bf x}_s) (Z^{{\bf x}}_s)^{-2}
\right) ds  = 0.
\end{equation}
For all sufficiently large $n$, the integrand in the last integral is estimated from below by $(\bar{\beta} - \bar{\alpha})/2$.
Therefore,
\[
\mathrm{P} ( \mathcal{{X}}^{\bf x}_{\sigma^{\bf x}} \in \Gamma_{n+1}) \geq \frac{1}{2} + 
\frac{c' (\bar{\beta} - \bar{\alpha})}{4},~~~
{\bf x} \in \Gamma_n,
\]
as required.
\qed
\\

In the proof of Lemma~\ref{sttm}, we saw that (\ref{slnn}) holds for sufficiently large $n$ if $\bar{\alpha} < \bar{\beta} $. The condition that $n$ is large was needed to ensure that the integrand in (\ref{innn}) was positive. If we momentarily consider $\rho = 0$ and observe that
the process $\Phi(\mathcal{{X}}^{\bf x}_t)$ is translation-invariant in the second variable in this case, then we obtain that
(\ref{slnn}) holds for all $n$. 

Now observe that the process $X^x_t$ in $(y,z)$ coordinates (i.e., without applying the transformation $\Psi_\varepsilon$) 
is governed by the same operator, up to the correction term $R$ (see formula (\ref{cott})), as the process $\mathcal{{X}}^{\bf x}_t$ with $\rho = 0$ 
in $(y, {\bf z})$ coordinates. For $n \in \mathbb{Z}$, define $\gamma_n  = \{(y,z): \psi(y) + \ln(z) = n\}$ - these are the analogues of
the sets $\Gamma_n$, but in $(y,z)$ coordinates. For For $x \in \gamma_n$, let 
\[
\tilde{\sigma}^{{x}} = \inf\{t \geq 0: {{X}}^{ x}_t \in \gamma_{n-1} \bigcup \gamma_{n+1}\}.
\]
Using the smallness of the coefficients of $R$ for small $z$,  it is easy to show, similarly to the proof of (\ref{slnn}), that, for all $n$ sufficiently close to $-\infty$, 
\[
\mathrm{P} ({{X}}^{ x}_{\tilde{\sigma}^{{x}} }  \in \gamma_{n+1}) \geq \frac{1}{2} + c,~~~
{ x} \in \gamma_n,
\]
provided that $\bar{\alpha} < \bar{\beta} $. Similarly, for $\bar{\alpha} < \bar{\beta} $,
\[
\mathrm{P} ({{X}}^{ x}_{\tilde{\sigma}^{{x}} }  \in \gamma_{n-1}) \geq \frac{1}{2} + c,~~~
{ x} \in \gamma_n.
\]
Since the process $X^x_t$ does not degenerate in $D$, these two inequalities immediately  imply the following lemma, which is not directly used in the proof of Theorem~\ref{mte}, but may be of independent interest.
\begin{lemma} \label{mlem}
If $\bar{\alpha} > \bar{\beta} $, then $\mathrm{P} (\lim_{t \rightarrow \infty} {\rm dist}(X^x_t, S) = 0) = 1$ for each $x \in \overline{D}$. 
If $\bar{\alpha} < \bar{\beta}$, then $\mathrm{P} (\lim_{t \rightarrow \infty} {\rm dist}(X^x_t, S) = 0) = 0$ for each $x \in {D}$. 
\end{lemma}

\vspace{0.8cm}

\noindent
{\it Proof of Lemma~\ref{mleee}.} If $S$ is attracting or neutral, we define $u({\bf x}) = \mathrm{E} f(\mathcal{X}^{\bf x}_{{\tau}^{\bf x}})$,
where the right-hand side is correctly defined since $\mathrm{P}({\tau}^{\bf x} < \infty) = 1$ (Lemma~\ref{sttm}). This is a standard
probabilistic representation of the solution to equation (\ref{fieq}), and the solution is unique in~$\mathcal{C}$ (see, e.g., \cite{F85}).

If $S$ is repelling, we define $h({\bf x}) = \mathrm{P} ({{\tau}^{\bf x}} < \infty)$. By Lemma~\ref{sttm},
$\lim_{{\bf z} \rightarrow \infty} \mathrm{P}({\tau}^{(y,{\bf z})} < \infty) = 0$, and thus $h \in  \mathcal{C}^2$ is 
the unique solution to (\ref{fieq2}).  Moreover, the process ${\widehat{\mathcal{X}}}^{\bf x}_t$, defined by conditioning 
${\mathcal{X}}^{\bf x}_t$ on the event 
$\{{{\tau}^{\bf x}} < \infty\}$, is governed by the operator $\widehat{M} u = h^{-1} M ( h u)$. The operator is non-degenerate, 
and, by construction, the process reaches $S$ with probability one for each initial point ${\bf x}$ (see, e.g., \cite{Pin}). Therefore, there is a unique solution
$u \in \mathcal{C}$ to (\ref{fieq3}), which is given by $u({\bf x}) = \mathrm{E} f(\widehat{\mathcal{X}}^{\bf x}_{\widehat{\tau}^{\bf x}})$,  where
$ \widehat{\tau}^{\bf x} = \inf\{t \geq 0: \widehat{\mathcal{X}}^{\bf x}_t \in S \times \{0\}\} $.

It remains to prove part (c) of the lemma. We will need the following fact. Suppose that the generator 
of a diffusion process $H^x_t$ is a uniformly elliptic operator in a bounded domain $G$ with a smooth boundary. Let $K \subset G$ 
be compact. Let $\tau^x = \inf\{t \geq 0: H^x_t \in \partial G$. Let   $\mu^x$ be the measure on $\partial G$ induced by $H^x_{\tau^x}$,
and let $p^x$ be its density with respect to the Lebesgue measure. Then there is a constant $c > 0$ such that
\begin{equation} \label{esbl}
p^x(\tilde{x}) \geq c,~~x \in K, \tilde{x} \in \partial G.
\end{equation}
The bound $c$ can be chosen to be the same for all the operators that have the same ellipticity constant and bound on the $C$-norm of the coefficients.  

First, consider the case when $S$ is attracting or neutral. 
Recall the definition of the sets $\Gamma_n$ and $R_n^+$ from (\ref{subsets}). For ${\bf x} \in R_n^+$, let
\[
\sigma_n^{{\bf x}} = \inf\{t \geq 0: \mathcal{{X}}^{\bf x}_t \in \Gamma_{n}\}.
\]
Then, for ${\bf x} \in R_n^+$,
\begin{equation} \label{fkf}
u({\bf x}) = \mathrm{E} u( \mathcal{{X}}^{\bf x}_{\sigma_n^{{\bf x}} } ).
\end{equation}
Let $V_n = \sup_{{\bf x}_1,{\bf x}_2 \in \Gamma_n}|u({\bf x}_1) - u({\bf x}_2)|$. Then $V_0 \leq  \sup_{x_1,x_2 \in S}|f(x_1) - f(x_2)|$ and,
by (\ref{fkf}), 
\[
\sup_{{\bf x}_1,{\bf x}_2 \in R^+_n}|u({\bf x}_1) - u({\bf x}_2)| \leq V_n.
\]
Thus it is sufficient to show that $V_n \rightarrow 0$ as $n \rightarrow \infty$. Since the operator $\mathcal{A}$ in (\ref{genaa}) is uniformly elliptic and
its coefficients are bounded (uniformly in $n \geq 0$) in the domain bounded by $\Phi(\Gamma_{n})$ and $\Phi(\Gamma_{n+2})$, (\ref{esbl}) 
is applicable to the process $\Phi(\mathcal{{X}}^{\bf x}_t)$  with  $K = \Phi(\Gamma_{n+1})$. Consequently, for the density 
$p^{\bf x}_n$ of the measure $\mu^{\bf x}_n$ induced by 
$ \Phi(\mathcal{{X}}^{\bf x}_{\sigma_n^{{\bf x}} })$ on $\Phi(\Gamma_n)$, we have 
$p^{\bf x}_n(\tilde{{\bf  x}}) \geq c$, 
${\bf x} \in \Gamma_{n+1}$, $\tilde{{\bf x}} \in  \Phi(\Gamma_n)$. (Here we used the fact that $p^{\bf x}_n \geq \overline{p}^{\bf x}_n$, where
$\overline{p}^{\bf x}_n$ corresponds to stopping the process $\mathcal{{X}}^{\bf x}_t$ 
on $\Gamma_n \bigcup \Gamma_{n+2}$ rather than on $\Gamma_n$.)
Therefore, by (\ref{fkf}), for ${\bf x}_1,{\bf x}_2 \in \Gamma_{n+1}$,
\[
|u({\bf x}_1)  -u({\bf x}_2)| \leq (1 - c \lambda(\Phi(\Gamma_n))) V_n,
\]
where $\lambda(\Phi(\Gamma_n)) = \lambda(\Phi(\Gamma_0))$ is the Lebesgue measure of $\Phi(\Gamma_n)$. Thus, $V_n \leq V_0 (1 - c 
\lambda(\Phi(\Gamma_0)))^n \rightarrow 0$ as $n \rightarrow \infty$, which implies that there is a limit
\[
\overline{u} = \lim_{{\bf z} \rightarrow \infty} u (y,{\bf z}) 
\]
uniformly in $y \in S$. 

The same argument applies in the case  when $S$ is repelling. We only need to observe that the operator governing the 
conditioned process $ \Phi(\widehat{\mathcal{X}}^{\bf x}_t)$ 
in the domain  between $\Phi(\Gamma_n)$ and $\Phi(\Gamma_{n+2})$ has coefficients bounded from above and the ellipticity constant bounded from below uniformly in $n \geq 0$,  as follows from the standard 
elliptic estimates on the function $h$. \qed
\\

In the proof of Lemma~\ref{mleee}, we saw that $u({\bf x}) = \mathrm{E} f(\mathcal{X}^{\bf x}_{{\tau}^{\bf x}})$ if $S$ is
attracting or neutral, and $u({\bf x}) = \mathrm{E} f(\widehat{\mathcal{X}}^{\bf x}_{\widehat{\tau}^{\bf x}})$ if $S$ is repelling. 
Let $\nu^{\bf x}$ be the measure on $S$ induced by $\mathcal{X}^{\bf x}_{{\tau}^{\bf x}}$ in the former case, or by 
$\widehat{\mathcal{X}}^{\bf x}_{\widehat{\tau}^{\bf x}}$ in the latter case. Thus $u({\bf x}) = \int_S f d \nu^{\bf x}$. Since the 
mapping $f \rightarrow \overline{u}$ is linear and continuous from $C(S)$ to $\mathbb{R}$, there is a measure $\nu$ on $S$ such that
$\overline{u} = \int_S f d \nu$. From part (c) of Lemma~\ref{mleee}, it follows that $\nu^{(y, {\bf z})} \rightarrow \nu$ weakly, as
${\bf z} \rightarrow \infty$. 

Similarly, assuming that Theorem~\ref{mte} holds, the measure on $S$ induced by $X^{x, \varepsilon}_{\tau^{x, \varepsilon}}$ converges
weakly to $\nu$ for each $x \in D$, as $\varepsilon \downarrow 0$, where $\tau^{x, \varepsilon}$ is the first time when $X^{x, \varepsilon}_t$
hits the boundary.

\section{Proof of the main result} \label{mresp}

In this section, we prove Theorem~\ref{mte}. Recall that $\mathcal{X}^{\bf x}_t$ is the process on $S \times [0,\infty)$ with the generator 
$M$ and $\mathcal{X}^{{\bf x},\varepsilon}_t = ({Y}^{{\bf x},\varepsilon}_t,{Z}^{{\bf x},\varepsilon}_t)$ is 
the process with the generator $M^\varepsilon$. For each $r > 0$, the latter operator can be defined on $S \times [0,r]$, 
provided that $\varepsilon > 0$  sufficiently small. 

 Let $K$ be a compact in $D$ and let $\delta > 0$. 
We will show that there is $\varepsilon > 0$ such that $|u^\varepsilon(x) - \overline{u}| \leq \delta$, $x \in K$,
provided that $\varepsilon \leq \varepsilon_0$. First, consider the case when  $S$ is attracting or neutral.

Recall the definition of $\Gamma_n$ from (\ref{subsets}).
Let $n \in \mathbb{N}$ be such that $|u({\bf x}) - \overline{u}| \leq \delta/2$ for
${\bf x} \in \Gamma_n$, where  $u({\bf x}) = \mathrm{E} f(\mathcal{X}^{\bf  x}_{{\tau}^{\bf x}})$
is the function from Lemma~\ref{mleee}.  Fix $\eta > 0$, to be specified later. Take $T, r > 0$ such that
\begin{equation} \label{bnn}
\mathrm{P}({\tau}^{\bf x} > T) \leq \eta,~~~ \mathrm{P}(\sup_{0 \leq t \leq {\tau}^{\bf x}} {Z}^{\bf x}_t > r-1)  \leq \eta,~~~
{\bf x}\in \Gamma_n.
\end{equation}
Such $T$ and $r$ exist  since $\mathrm{P}({\tau}^{\bf x} < \infty) = 1$ and the probabilities in the left-hand side of
both inequalities depend continuously on ${\bf x}$. Define
\[
{\tau}^{{\bf  x},\varepsilon}_r =  \begin{cases} \inf\{t \geq 0: {Z}^{{\bf x},\varepsilon}_t  = 0\}, & 
{\rm if }~ \inf\{t \geq 0: {Z}^{{\bf x},\varepsilon}_t  = 0\} < \inf\{t \geq 0: {Z}^{{\bf  x},\varepsilon}_t  = r\}, \\ \infty, & {\rm otherwise}. \end{cases}
\]
Since $M^\varepsilon$ is a small perturbation of $M$ on $S \times [0,r]$
(formula (\ref{genpp0})), from (\ref{bnn}) it follows that, for all sufficiently small $\varepsilon$,
\begin{equation} \label{bnn2}
\mathrm{P}\left({\tau}^{{\bf x}, \varepsilon}_r < r,~\|\mathcal{X}^{{\bf x}}_{{\tau}^{{\bf x}}} - 
\mathcal{X}^{{\bf x},\varepsilon}_{{\tau}^{{\bf x},\varepsilon}_r} \| \leq \eta \right) \geq 1-  2\eta,~~~{\bf x} \in \Gamma_n.
\end{equation}
Observe that
\[
|u^\varepsilon({\bf x}) - \mathrm{E} \left( f(\mathcal{X}^{{\bf x},\varepsilon}_{{\tau}^{{\bf x},\varepsilon}_r});~ 
{\tau}^{{\bf x}, \varepsilon}_r < r \right) | \leq 2 \eta \sup |f|,
\]
\[
|u({\bf x}) - \mathrm{E} \left(
f(\mathcal{X}^{{\bf x}}_{{\tau}^{{\bf x}}}),~{\tau}^{{\bf x}, \varepsilon}_r < r \right) |\leq 2 \eta \sup |f|. 
\]
(Here, the solution $u^\varepsilon$ to equation (\ref{direq}) is considered in $(y,{\bf z})$ coordinates.)
Therefore, from (\ref{bnn2}) it follows that, for all sufficiently small $\varepsilon$,
\[
 |u^\varepsilon({\bf x}) -  u({\bf x}) | \leq \mathrm{E} \left( |f(\mathcal{X}^{{\bf x},\varepsilon}_{{\tau}^{{\bf x},\varepsilon}_r}) - 
f(\mathcal{X}^{{\bf x}}_{{\tau}^{{\bf x}}}) |;~ 
{\tau}^{{\bf x}, \varepsilon}_r < r \right) | + 4 \eta \sup |f| \leq
\]
\[
\sup_{x_1, x_2 \in S, \|x_1 - x_2\| \leq \eta} |f(x_1)  -f(x_2)| +  8 \eta \sup |f| \leq \frac{\delta}{2},~~~{\bf x} \in \Gamma_n,
\]
where the last inequality is obtained by selecting a sufficiently small $\eta$. Thus,
\[
|u^\varepsilon({\bf x}) -  \overline{u}| \leq | u({\bf x}) - \overline{u}| +  |u^\varepsilon({\bf x}) -  u({\bf x}) |  \leq \delta, ~~~{\bf x} \in \Gamma_n.
\]

Next, consider the case when $S$ is repelling. Again, let $n \in \mathbb{N}$ be such that $|u({\bf x}) - \overline{u}| \leq \delta/2$ for
${\bf x} \in \Gamma_n$. Fix $\eta > 0$, to be specified later.  For $r > 0$, let
\[
{\tau}^{{\bf x}}_r =  \begin{cases} \inf\{t \geq 0: {Z}^{{\bf x}}_t  = 0\}, & 
{\rm if }~ \inf\{t \geq 0: {Z}^{{\bf x}}_t  = 0\} < \inf\{t \geq 0: {Z}^{{\bf x}}_t  = r\}, \\ \infty, & {\rm otherwise}. \end{cases}
\]
Let $r$ be sufficiently large so that 
\[
|u({\bf x}) - \mathrm{E} \left( f(\mathcal{X}^{{\bf x}}_{{\tau}^{{\bf x}}_r}) | {\tau}^{{\bf x}}_r < \infty \right)| \leq \eta
\]
for all ${\bf x} \in \Gamma_n$. From the proximity of $\mathcal{X}^{{\bf x}}_t$ and $\mathcal{X}^{{\bf x},\varepsilon}_t$, 
using the same arguments as above, it is easy to show that 
\[
|\mathrm{E} \left( f(\mathcal{X}^{{\bf x},\varepsilon}_{{\tau}^{{\bf x},\varepsilon}_r}) | {\tau}^{{\bf x},\varepsilon}_r < \infty \right) -  
\mathrm{E} \left( f(\mathcal{X}^{{\bf x}}_{{\tau}^{{\bf x}}_r}) | {\tau}^{{\bf x}}_r < \infty \right) | \leq \eta,~~~{\bf x} \in \Gamma_n, 
\]
for all sufficiently small $\varepsilon$, and therefore 
$|\mathrm{E} \left( f(\mathcal{X}^{{\bf x},\varepsilon}_{{\tau}^{{\bf x},\varepsilon}_r}) | {\tau}^{{\bf x},\varepsilon}_r < \infty \right)  - u({\bf x})| \leq 2 \eta$.  Representing 
$u^\varepsilon({\bf x}) = \mathrm{E} ( f(\mathcal{X}^{{\bf x},\varepsilon}_{{\tau}^{{\bf x},\varepsilon}}))$ as a sum of contributions
from successive visits to $\Gamma_n$ after reaching the surface defined by $\{{\bf z} = r\}$, we obtain 
$|u^\varepsilon({\bf x})   - u({\bf x})| \leq 2 \eta \leq \delta/2$, where the last inequality follows by taking a sufficiently small $\eta$. Thus,
$|u^\varepsilon({\bf x}) -  \overline{u}| \leq \delta$ for ${\bf x} \in \Gamma_n$. 

For $x \in K$, the estimate follows from the probabilistic representation of the solution $u^\varepsilon$ in the
domain bounded by $\Gamma_n$. \qed

\section{Remarks and generalizations} \label{gener}

Boundary problems for elliptic partial differential equations often appear as a result of stabilization, as time tends to infinity,
in initial-boundary value problems for the corresponding evolutionary equations. For instance, the Dirichlet problem 
\[
L^\varepsilon u^\varepsilon (x)  = 0,~~x \in D;~~~~u^\varepsilon(x)  = \psi(x),~~x \in \partial D,
\]
arises when one considers the first initial-boundary value problem
\begin{equation} \label{direq2}
\frac{\partial u^\varepsilon (t,x)}{\partial t} = L^\varepsilon u^\varepsilon (t, x),~~t > 0, x \in D;
\end{equation}
\[
u^\varepsilon(0, x)  = 
g(x),~~x \in D;~~~~u^\varepsilon(t,x) = \psi(x),~~t > 0, x \in \partial D.
\]
Here, we assume that $g \in C(\overline{D})$, $\psi \in C(\partial D)$.  The second initial-boundary value problem
\begin{equation} \label{direq3}
\frac{\partial  u^\varepsilon  (t,x)}{\partial t} = L^\varepsilon u^\varepsilon (t, x),~~t > 0, x \in D;
\end{equation}
\[
u^\varepsilon(0, x)  = 
g(x),~~x \in D;~~~~\frac{\partial u^\varepsilon(t,x)}{\partial n^\varepsilon(x)} = 0,~~t > 0, x \in \partial D,
\]
where  $n^\varepsilon(x)$  is the co-normal to $\partial D$ at $x$,
leads to the Neumann problem.  Here and below, we assume that $\partial D = S_1 \bigcup ... \bigcup S_m$, where $S_1,...,S_m$ are
disjoint smooth connected $(d-1)$-dimensional manifolds. 

If the operator depends on a small parameter $\varepsilon$, the limiting behavior of the solution of the initial-boundary
value problem as $t \rightarrow \infty$ and, simultaneously, $\varepsilon \downarrow 0$ should be considered. This double limit, in general, may not exist; the limiting behavior depends on how the point $(\varepsilon^{-1}, t(\varepsilon))$ approaches infinity. The solutions of initial-boundary value problems (\ref{direq2}) and (\ref{direq3}) can be written as expectations of certain functionals of the
corresponding diffusion processes that depend on the parameter $\varepsilon$, and the dependence of 
the limit of $u^\varepsilon(t(\varepsilon),x)$ on the asymptotics of $t(\varepsilon)$ is a
manifestation of metastability for the underlying diffusion (see, e.g., \cite{FW}).

Below, we briefly discuss the asymptotic behavior for solutions to parabolic equations in the case when the boundary has several connected 
components. This topic (and the proofs of the claims made below) will be the subject of a forthcoming paper. 

As before, we can classify each component of the boundary as attracting, neutral, or repelling, based on the behavior of the coefficients
of the operator near the boundary. Moreover, with each attracting component $S_k$, we can associate a number $\gamma_k >0$ such
that the time it takes the process $X^{x,\varepsilon}_t$ starting at $x \in S_k$ to exit a fixed neighborhood of $S_k$ is of order
$\varepsilon^{-\gamma_k}$ when $\varepsilon \downarrow 0$. If $S_k$ is repelling, we can associate a number $\tilde{\gamma}_k$ to it
such that the time it takes the process $X^{x,\varepsilon}_t$, $x \in D$, conditioned on exiting the domain through $S_k$,  to reach $S_k$ is of order $\varepsilon^{-\tilde{\gamma}_k}$.  
 The numbers $\gamma_k$ (if $S_k$ is attracting) and $\tilde{\gamma}_k$ (if
$S_k$ is repelling) can be
found by solving a certain spectral problem that involves the operator $L$ restricted to $S_k$ and the leading terms of the coefficients
of the operator near~$S_k$.

As suggested  in \cite{Finf}, the long-time behavior of a perturbed system  (of the diffusion process with a small parameter in our case) can be described  by a motion on the simplex of invariant probability measures of the unperturbed system. The vertices of this simplex are the ergodic probability measures of the unperturbed system. The limiting behavior of $X^{x,\varepsilon}_{t(\varepsilon)}$ as $\varepsilon \downarrow 0$ and $t(\varepsilon) \rightarrow \infty$ can be described by the time evolution of two ``coordinates". The first coordinate
is a point $\mu_{X^{x,\varepsilon}_{t(\varepsilon)}}$ of the simplex $M$ of invariant probability measures of the non-perturbed system, 
where $\mu_y$, $y \in \overline{D}$, denotes the limiting distribution of the non-perturbed system starting at $y$. In our case, such a limit exists  for each $y \in \overline{D}$. The second coordinate is a point $r_{X^{x,\varepsilon}_{t(\varepsilon)}}$ in the support of 
$\mu_{X^{x,\varepsilon}_{t(\varepsilon)}}$ that is nearest to $X^{x,\varepsilon}_{t(\varepsilon)}$. Under appropriate assumptions on 
the time scale $t(\varepsilon)$, the measures $\mu_{X^{x,\varepsilon}_{t(\varepsilon)}}$  have a limit $\overline{\mu}^x$ (that depends on the initial
point and the time scale). The second coordinate, $r_{X^{x,\varepsilon}_{t(\varepsilon)}}$, converges in distribution to the measure 
$\overline{\mu}^x$  as $\varepsilon \downarrow 0$. If the set of stable (in a certain sense) invariant probability measures for the 
unperturbed system is finite, then  these stable measures serve as limits for $\mu_{X^{x,\varepsilon}_{t(\varepsilon)}}$ at different time
scales and
the transitions between them occur as jumps. In the case of parabolic problems being discussed, these switches between different
metastable states for the underlying diffusion are manifested by the fact that, for each $x \in D$, the
limit $\lim_{\varepsilon \downarrow 0} u^\varepsilon(t(\varepsilon), x)$ can take values from a finite set if
we disallow certain ``transitional" time scales.

In the case of the first initial-boundary value problem, the
process $\overline{X}^{x, \varepsilon}_t$ in $\overline{D}$, obtained from ${X}^{x, \varepsilon}_t$ by  stopping it when 
it hits the boundary,  should be considered:
\[
\overline{X}^{x, \varepsilon}_t =  \begin{cases} {X}^{x, \varepsilon}_t, & 
t \leq \tau^{x,\varepsilon}, \\ {X}^{x, \varepsilon}_{ \tau^{x,\varepsilon}}, & t > \tau^{x,\varepsilon}, \end{cases}
\]
where $\tau^{x,\varepsilon} = \inf\{t \geq 0: {X}^{x, \varepsilon}_t \in \partial D \}$. Note that $\mathrm{P}( \tau^{x,\varepsilon} 
< \infty) = 1$ for each $\varepsilon > 0$, $x \in \overline{D}$. The unperturbed process in this case is 
\[
\overline{X}^{x}_t =  \begin{cases} {X}^{x}_t, & 
{\rm if}~ x \in D, \\ x, & {\rm if}~ x \in \partial D. \end{cases}
\]
Each $\delta$-measure concentrated at a point $x \in \partial D$ is an ergodic invariant probability measure for $\overline{X}^{x}_t$,
and the simplex of all invariant probability measures coincides with with the set of all probability measures on $\partial D$ if at least 
one component of $\partial D$ is attracting. If all the components of the boundary are repelling, then $\overline{X}^{x}_t$ has, in addition,
an invariant probability measure $\mu$ that is absolutely continuous with respect to the Lebesgue measure  on $D$. For each attracting 
component $S_k \subseteq \partial D$, one can define a probability measure $\nu_k$ similarly to the measure $\nu$ 
in Section~\ref{homog} (it is the limiting distribution, as $\varepsilon \downarrow 0$, of the exit point for the process
${X}^{x, \varepsilon}_t$ conditioned on exiting the domain through $S_k$). We can also define $\pi_k$ as the
invariant measure of the process $X^x_t$ restricted to $S_k$.

It turns out 
that if at least one $S_k$ is attracting and $1 \ll t(\varepsilon) \ll |\ln(\varepsilon)|$, then  
$\lim_{\varepsilon \downarrow 0} u^\varepsilon(t(\varepsilon), x)$,   exists and \be 
is equal to a certain linear combination of  $\int_{S_k} g  d \pi_k$ over all attracting $S_k$. The coefficients in this linear combination depend on $x \in D$ and can be calculated as solutions of the appropriate Dirichlet problems for the
operator $L$. The condition $t(\varepsilon) \ll |\ln(\varepsilon)|$ appears because it takes time of order $|\ln(\varepsilon)|$ for the process $X^{x,\varepsilon}_t$ to reach  the boundary if there is at least one attracting component. 

If at least one $S_k$ is attracting and $t(\varepsilon) \gg |\ln(\varepsilon)|$, then 
$\lim_{\varepsilon \downarrow 0} u^\varepsilon(t(\varepsilon), x)$ is equal to the linear combination (with the same coefficients as in the case above) of the quantities $\int_{S_k} \psi  d \nu_k$. Thus the asymptotic 
behavior $u^\varepsilon(t(\varepsilon), x)$  switches at times of order $|\ln(\varepsilon)|$.

If all $S_k$ are repelling, a switch also happens, but at a different time scale. Let $k^*$ be such that $\tilde{\gamma}_{k^*} =
\min_k \tilde{\gamma}_k$. For simplicity, we assume that such $k^*$ is unique. Then 
$\lim_{\varepsilon \downarrow 0} u^\varepsilon(t(\varepsilon), x) = \int_D g d \mu$ if 
$1 \ll t(\varepsilon) \ll \varepsilon^{-\tilde{\gamma}_{k^*}}$, and 
$\lim_{\varepsilon \downarrow 0} u^\varepsilon(t(\varepsilon), x) = \int_{S_{k^*}} \psi d \nu_{k^*}$ if 
$t(\varepsilon) \gg \varepsilon^{-\tilde{\gamma}_{k^*}}$.

Now let us discuss the second initial-boundary value problem (\ref{direq3}). In this case, the process corresponding to the perturbed problem is the diffusion $\hat{X}^{x,\varepsilon}_t$ in $\overline{D}$ governed by the operator $L^\varepsilon$ inside $D$ and reflecting instantaneously in the
co-normal direction on $\partial D$. In this case, the unperturbed process is the diffusion $X^x_t$, defined in Section~\ref{intro}, 
restricted to $\overline{D}$. One ergodic probability measure $\pi_k$ of this process is concentrated on each connected component of
the boundary $S_k \subseteq \partial D$. If at least one component of the boundary is attracting, then the simplex of invariant measures
is the convex envelope of the measures $\pi_k$, $1 \leq k \leq m$. If all $S_k$ are repelling, then the abosolutely
continuous measure $\mu$ on $D$ (which was introduced above) serves as an additional invariant probability measure for the process $X^x_t$ in $\overline{D}$. The
simplex of invariant probability measures in this case is the convex
envelope of the collection $\mu, \pi_1,...,\pi_m$. 

If there are attracting components of the boundary, then we denote them by $S_1,...,S_m$, arranged in such an order that 
$\gamma_1 < \gamma_2 < ... < \gamma_m$ (for simplicity, we assume that all $\gamma_k$ are distinct). We also put $\gamma_0 = 0$. For $\varepsilon^{-{\gamma}_{k-1}} \ll t(\varepsilon) \ll \varepsilon^{-{\gamma}_{k}}$, $1 \leq k \leq 
m$, the limit $\lim_{\varepsilon \downarrow 0} u^\varepsilon(t(\varepsilon), x)$ can be expressed as a linear combination of the quantities
$\int_{S_i} g d \pi_i$ that extends over $i$ satisfying $k \leq i \leq m$. The coefficients in the linear combination depend on $x \in D$ and
can be found by solving the corresponding auxiliary elliptic problems. For $t(\varepsilon) \gg \varepsilon^{-{\gamma}_{m}}$,
$\lim_{\varepsilon \downarrow 0} u^\varepsilon(t(\varepsilon), x) = \int_{S_m} g d \pi_m $, which does not depend on $x$. Here, as we see,
the switch in the asymptotic behavior of $u^\varepsilon(t(\varepsilon), x) $ happens at several time scales.
 
Finally, if all the components of the boundary are repelling and $t(\varepsilon) \gg 1$, then $\lim_{\varepsilon \downarrow 0} u^\varepsilon(t(\varepsilon), x) = \int_D g d \mu$.
\\

Throughout this paper, it was assumed 
that the restriction of the operator $L$ to the boundary was a non-degenerate diffusion. Actually, what is
really important is that the the diffusion process on $S$ (or each component of $S$ in the case of non-connected boundary) has a unique
invariant probability measure. For instance, it is sufficient to assume that the process satisfies the Doeblin condition. Multiplicity of
the invariant probability measures leads to more sophisticated limiting behavior of $u^\varepsilon(t(\varepsilon), x)$.

One can also consider the case when the operator $L$ has degeneracies inside the domain $D$. Again, the structure of the simplex of invariant
probability measures of the non-perturbed process will play an important role in describing  the limiting behavior of 
$u^\varepsilon(t(\varepsilon), x)$ as $\varepsilon \downarrow 0$, $t(\varepsilon) \rightarrow \infty$.
\\
\\
\noindent {\bf \large Acknowledgments}: The work of  L. Koralov was supported by the Simons Foundation Fellowship (award number 678928)
and by the Russian Science Foundation,  project ${\rm N}^o$~20-11-20119. 
\\
\\

\end{document}